\documentstyle{amsppt}
\expandafter\let\csname logo\string @\endcsname=\empty
\magnification=1200
\pagewidth{5.4in}
\pageheight{7.5in}
\NoRunningHeads
\NoBlackBoxes
\define\a{\alpha}
\predefine\barunder{\b}
\redefine\b{\beta}

\define\cd{\cdot}
\predefine\dotunder{\d}
\redefine\d{\delta}
\predefine\dotaccent{\D}
\redefine\D{\Delta}
\define\e{\epsilon}

\define\f{\frac}
\define\g{\gamma}
\define\G{\Gamma}

\define\lb\{{\left\{}
\define\la{\lambda}

\define\lra{\longrightarrow}

\define\oper{\operatorname}

\define\p{\partial}
\define\rb\}{\right\}}
\define\s{\sigma}

\define\sub{\subheading}

\define\vp{\varphi}

\define\x{\times}

\define\({\left(}
\define\){\right)}
\define\[{\left[}
\define\]{\right]}
\define\<{\left<}
\define\>{\right>}
\def\slantline#1#2#3#4#5{\hbox to 0pt{
}}


\def\SH{{\Cal H}}

\def\SL{{\Cal L}}
\def\SM{{\Cal M}}

\def\SO{{\Cal O}}

\def\SQ{{\Cal Q}}
\def\SR{{\Cal R}}
\def\SS{{\Cal S}}
\def\ST{{\Cal T}}


\def\BC{{\bold C}}

\def\BH{{\bold H}}

\def\BR{{\bold R}}


\def\bbc{{\Bbb C}}

\def\bbh{{\Bbb H}}

\def\bbr{{\Bbb R}}
\def\bbs{{\Bbb S}}

\def\bbz{{\Bbb Z}}

\def\iny{\exp(2\pi iny/\ell)}
\def\copi{\cosh\frac{\pi ns}{\ell}}
\def\sipi{\sinh\frac{\pi ns}{\ell}}
\def\gr{\text{gr}}
\def\grg{\text{Gr}_{s\g}}
\def\ungrg{\text{UnGr}_{s\g}}
\def\gro{\text{gr}(\s_0)}
\def\grt{\text{gr}(\s_t)}
\def\id{\oper{id}}
\def\tec{Teichm\"uller\ }
\def\elg{\ell_\g(\text{gr}(\s_t))}
\def\ddt{\frac d{dt}\bigm|_{t=0}}
\def\ds{ds_{\text{gr}(\s_0)}}
\def\dst{ds_{\text{gr}(\s_t)}}
\def\barcn{\bar c_n}
\def\bardn{\bar d_n}
\def\dcyl{\frac{s}{2}}
\def\bigo{\SO}
\def\desitter{\bbs^3_1}
\def\complex{\bbc}

\def\rem{\medskip\noindent{{\bf Remark:\ \ }}}

\def\prf{\medskip\noindent{{\bf Proof:\ \ }}}
\def\pf#1{\medskip\noindent{{\bf Proof of #1.\ \ }}}
\def\diffo{\oper{Diff}_o}
\def\im{\oper{Im}}
\def\ree{\oper{Re}}

\topmatter
\title The grafting map of Teichm\"uller space
\endtitle
\affil  Kevin P. Scannell\\
Dept\. of Mathematics and Computer Science\\
St\. Louis University\\
St\. Louis, MO 63103\\
\\
Michael Wolf*\\
Dept\. of Mathematics\\
Rice University\\
Houston, TX 77251 
\endaffil
\thanks *Partially supported by NSF Grants DMS 9626565, DMS 9707770 (SCREMS)
\endthanks
\endtopmatter

\document

\centerline{October 11, 1998}

\sub{\S1. Introduction}

\sub{1.1 Statement and Context}

One of the underlying principles in the study of Kleinian groups is that
aspects of the complex projective geometry of quotients of $\hat\BC$
by the groups reflect properties of the three-dimensional hyperbolic 
geometry of the quotients of $\BH^3$ by the groups.
Yet, even though
it has been over thirty-five years since
Lipman Bers wrote down a holomorphic embedding of the \tec
space of Riemann surfaces in terms of the projective
geometry of a \tec space of quasi-Fuchsian manifolds, no 
corresponding parametrization in terms of the three-dimensional
hyperbolic geometry has been presented.  One of the goals of this paper is to
give such a parametrization. This parametrization is straightforward
and has been expected for some time (\cite{Ta97}, \cite{Mc98}): 
to each member of a
Bers slice of the space $QF$ of quasi-Fuchsian 3-manifolds, we 
associate the bending measured lamination of the convex hull facing the 
fixed ``conformal'' end.

The geometric relationship between a boundary component of a convex hull 
and the projective surface at infinity for its end is given by 
a process known as {\it grafting}, an operation on projective
structures on surfaces that traces its roots back at least to Klein
\cite{Kl33;\S50, p. 230}, with a modern history
developed by many authors (\cite{Ma69}, \cite{He75}, 
\cite{Fa83}, \cite{ST83}, \cite{Go87},\cite{GKM95},\cite{Ta97},\cite{Mc98}).
The main technical tool in our proof that bending measures give
coordinates for Bers slices, and the 
second major goal of this paper, is the  completion of the
proof of the ``Grafting Conjecture''. This conjecture
states that for a fixed measured
lamination $\lambda$, the self-map of \tec space induced by grafting 
a surface along $\lambda$ is a homeomorphism of \tec space; our
contribution to this argument is a proof of the injectivity 
of the grafting map. While the principal application of this 
result that we give is to geometric coordinates on the Bers slice
of $QF$, one expects that 
the grafting homeomorphism
might lead to other systems of 
geometric coordinates for other families of Kleinian groups
(see \S5.2);
thus we feel that this result is of interest in its own right.

A difficulty in proving injectivity results for maps of 
\tec space is that \tec space is a quotient space with 
no canonical sections; our approach is to choose a section 
in the space of metrics over \tec space defined via harmonic
maps.  Indeed,
our principal tool in proving the grafting conjecture is a study
of the the differential equation governing
the infinitesimal form of the energy density of such maps; this
study is complicated somewhat by the grafted metrics having a
mild singularity, and the infinitesimal form having a more serious
singularity.  Nevertheless, this equation is amenable to 
nearly a complete solution, and it is estimates based on this solution
which are the technical linchpins of our argument.

We now state our results and methods more precisely.
Throughout, $S$ will denote a fixed differentiable surface
which is closed, orientable, and of genus $g \geq 2$.    Let $T_g$
be the corresponding Teichm\"uller space of marked conformal 
structures on $S$, and let $P_g$ denote the deformation
space of (complex) projective structures on $S$ (see \S2 for
definitions).

There are two well-known parametrization of
$P_g$, each reflecting a 
different aspect of the general theory.
The first uses the Schwarzian derivative of the developing map
to obtain a quadratic differential on $S$,
holomorphic with respect to the complex structure underlying
the given projective structure.   This identifies $P_g$ with
the total space of the bundle $Q_g \to T_g$ of holomorphic
quadratic differentials over Teichm\"uller space.
This identification is representative of the complex analytic
side of the theory; see for instance
\cite{Ea81}, \cite{Gu81}, \cite{He75}, \cite{Kr69}, \cite{Kr71}, 
\cite{KM81}, \cite{MV94}, \cite{Sh87}, \cite{ST95}.

The second parametrization is due to Thurston and is more 
geometric in nature.  To describe it, fix a hyperbolic metric
$\sigma \in T_g$, a simple closed geodesic $\gamma \subset S$
of length $\ell$, and a positive real number $s$. 
Let 
$$
\tilde{A_s} = \{ (r,\theta) \in \bbc^* : | \theta - \pi / 2 | \leq s \} 
$$
and 
$$
A_s = \tilde{A_s} / \langle z \mapsto e^\ell z \rangle.
$$
Of course, if $s \geq 2 \pi$, we must interpret 
the projective structure on $\tilde{A_s}$ as being defined 
by a developing map which is no longer an embedding;
in any case, we call $A_s$ a (projective) {\it $s$-annulus}. 
A new projective structure on $S$ is defined by cutting the original
hyperbolic surface $(S,\sigma)$ open along 
the simple closed curve $\gamma$ and 
gluing in $A_s$.    This is the
{\it grafting} operation; 
it provided the first examples \cite{Ma68} of projective structures
for which the developing map is not a covering of its image.
Grafting extends by continuity from pairs $(\gamma, s)$
to general measured laminations, defining a map
$\Theta : \Cal{ML} \times T_g \to P_g$.   Thurston has shown
(in unpublished work) that $\Theta$ is a homeomorphism 
(see \cite{KT92}, \cite{La92}).

A natural problem is to understand how these geometric and
complex analytic aspects interact.  For instance, a
measured lamination $\lambda \in \Cal{ML}$ 
defines a slice $\Theta(\{ \lambda \} \times T_g) \subset P_g$; 
following this inclusion with the projection $P_g \to T_g$
defines a self-map of Teichm\"uller space
$\text{Gr}_\lambda : T_g \to T_g$.
Our main result can be stated concisely as follows:

\proclaim{Theorem A} $\text{Gr}_\lambda$ is a homeomorphism.
\endproclaim

This result was obtained in special cases
by McMullen \cite{Mc98} (one-dimensional Teichm\"uller spaces), 
and Tanigawa \cite{Ta97} (for integral points of $\Cal{ML}$,
using a result of Faltings \cite{Fa83}). Our result will hold for
all elements of $\SM\SL$ and all \tec spaces of 
finitely punctured Riemann surfaces of finite genus. (For the
sake of expositional ease, we write the proof for \tec spaces of
closed Riemann surfaces, but the extension to \tec spaces of
finitely punctured surfaces is mostly a matter of additional
notation: see the remark at the end of \S4.)

Theorem A allows one to understand various complex analytic
constructions in the theory of Teichm\"uller spaces and
Kleinian groups in terms of measured geodesic laminations
and the grafting construction.    As an example, 
we obtain the following corollary in \S 5.1:

\proclaim{Corollary} Let $B_Y$ be a Bers slice with fixed conformal
structure $Y$, and define a map $\beta : B_Y \to \Cal{ML}$ which assigns
the bending lamination on the component of the convex hull boundary
facing $Y$.  Then $\beta$ is a homeomorphism onto its image.
\endproclaim

The space of projective structures is
intimately related with the space of locally convex pleated maps
of $\tilde{S}$ into $\bbh^3$ (as detailed for instance in \cite{EM87}).
The dual notions are explored in \cite{Sc96}, 
where it is shown that $P_g$ classifies causally trivial
de~Sitter structures on $S \times \bbr$; here the grafting operation
corresponds to a ``stretching'' of the causal horizon.
We give an application of Theorem A to this situation in \S 5.3.

Finally, in \cite{Mc98} McMullen observes that Theorem A follows 
from the conjectural rigidity of hyperbolic cone $3$-manifolds 
(see \cite{HK98}); hence our result can be viewed as further positive 
evidence for the validity of this conjecture.

\sub{1.2 Outline of the Argument}
It has been shown that $\text{Gr}_\lambda$ is real analytic \cite{Mc98}
and proper \cite{Ta97}, therefore it suffices, as $T_g$ is a cell,
to prove local
injectivity.   For simplicity, assume the measured lamination $\lambda$
is given by some simple closed curve $\gamma$ on $S$ and
a non-zero transverse measure $s \in \BR_+$.
Fix a hyperbolic metric $\s_0 \in T_g$
and a small deformation $\s_t$  of $\s_0$ with the property
that the grafted surfaces are conformal; 
i.e. $\text{Gr}_\lambda(\s_0) = \text{Gr}_\lambda(\s_t)$ in $T_g$.
Each grafted surface may be equipped with a $C^{1,1}$ metric
$\grt$ which is flat on the inserted cylinder and hyperbolic elsewhere
(see \S2.2).    Thus, it is a consequence of the singular 
harmonic maps theory of \cite{GS92} that 
for every $t$ there is a unique harmonic
map $w_t : (S, \grt) \to (S, \gro)$ homotopic to the
identity; indeed we may adjust the metrics
by an isotopy and assume that $w_t$ is the identity map 
(it is straightforward that $w_t$ is a homeomorphism).
After this normalization, the two boundary curves
of the inserted cylinder may move around $S$ as $t$ varies;
the variation vector fields $V^\ell$ and $V^r$
so-defined (one for each of the two boundary components
$\gamma^\ell$ and $\gamma^r$ of the inserted cylinder)
are key pieces of information in the proof.

With this setup, the proof proceeds by 
first recognizing that the conformal factor 
$\SH_t = \frac{gr(\s_0)}{gr(\s_t)}$ has a second role as the 
holomorphic energy density of the harmonic map $w_t$; in that
second role, it satisfies the Bochner equation 
(given as equation (3.1.3) below).   
It turns out to be easier to analyze the linearized equation for $\dot \SH$ 
at $t = 0$ (equation (3.1.5)).

Our strategy for solving (3.1.5) is straightforward.
We think of (3.1.5) as representing two different equations;
the first on the open inserted flat cylinder $S_0$ (where $K = 0$)
and the second on the cut-open hyperbolic surface $S_{-1}$
(where $K = -1$).
A general solution to the first equation can be found easily.
To study the solution to the second equation and the 
global solution to (3.1.5) (note that the complete equation
(3.1.5) contains a term given as a measure $2\dot K$ supported
on $\gamma^\ell$ and $\gamma^r$), we begin with an observation:
the normal derivative of $\dot \SH$ across the boundary curves 
$\gamma^\ell$ and $\gamma^r$ appears
as the inhomogeneous term 
in an ordinary differential equation (see \S2.4) 
which can be solved for the variational fields $V^{\ell}$ and $V^r$.
These vector fields $V^{\ell}$ and $V^r$ in turn determine the 
normal derivatives of the global solution $\dot \SH$
as computed from the hyperbolic side $S_{-1}$.
Finally, integrating $\dot \SH \Delta \dot \SH$ by parts on $S_{-1}$
and using our knowledge of the boundary terms forces 
$\dot \SH$ to vanish identically. 
This then implies that the original metrics, $\s_t$ and $\s_0$,
are infinitesimally conformal, proving the desired 
local injectivity.

It is helpful in understanding the overall argument to 
note that we ignore the fact that the length of
the inserted cylinder (whose length is always denoted $s$) is constant in $t$
until the very end of the proof.  
This is discussed in some
detail in \S 3.3.

The case where $\lambda$ is a general measured lamination and not
just a simple closed curve (or a system of disjoint simple closed
curves) follows from approximating the general lamination by simple
closed curves, approximating the conformal deformation $Gr_{\lambda}(\s_t)$ by
quasi-conformal deformations $Gr_{s_m \g_m}(\s_t)$,
and then extending the previous argument for the simple closed
curves and conformal deformations 
to find identities involving only quantities that are continuous
on the space $\SM \SL$.

\sub{Acknowledgments} The authors wish to acknowledge several
pleasant and very useful conversations with
John Polking on regularity issues, with Robert Hardt on
properties of harmonic maps to singular spaces, and with 
Jim Anderson on his and Dick Canary's work on limits of
Kleinian groups.

\sub{\S2. Notation and Background}

\sub{2.1 Teichm\"uller Space, Bers embedding} Let $S$ denote a
smooth surface of genus $g\ge2$, and let $\SM_{-1}=\SM_{-1}(S)$ denote
the space of metrics $\rho|dw|^2$ on $S$ with Gaussian curvature
identically $-1$. The group $\diffo$ of diffeomorphisms of $S$ homotopic
to the identity acts on $\SM_{-1}$ by pullback: if $\phi\in\diffo$, then
$\phi\cd\rho=\phi^*\rho$. We define the Teichm\"uller space of genus
$g$, $T_g$, to be the quotient space $T_g=\SM_{-1}/\diffo$, i.e.,
equivalence classes of metrics in $\SM_{-1}$ under the action of
$\diffo$. A metric $(S,\rho)$ represents a conformal class of metrics
on S, hence a Teichm\"uller equivalence class of Riemann surfaces.
Let $QD(\s)$ denote the $3g-3$ dimensional complex vector space of
holomorphic quadratic differentials on $(S,\rho)$.

There are a number of continuous and real-analytic parametrizations
of the Teichm\"uller space $T_g$ and one complex analytic 
parametrization given by Lipman Bers \cite{Be64}.  The Bers embedding, as it 
is usually known (see \cite{Na88} for a comprehensive account), is 
defined as follows.  Fix a point $Y$ in $T_g$.  Then, for any
(variable) point $X\in T_g$, consider the quasi-Fuchsian
manifold $Q(X,Y)$ with conformal boundaries $X$ and $Y$
and fundamental
group $\G(X,Y)$. There is a 
simultaneous uniformization homeomorphism $F:\hat\BC \to \hat\BC$
of the sphere $\hat\BC$ 
which does the following: 1) it equivariantly and conformally
maps the unit disk $\Delta$ to the universal cover of $Y$,
2) it equivariantly and quasi-conformally 
maps the complement $\Delta^*$ of the unit disk to the universal cover
of $X$, and 3) it conjugates $\G(Y,Y)$ to $\G(X,Y)$. As
$F\Bigm|_{\Delta}$ is conformal, we may take its Schwarzian
derivative, say $\SS(F\Bigm|_{\Delta})= \Psi_X$. The holomorphic 
function $\Psi_X$ descends to a holomorphic quadratic differential 
on the Riemann surface $Y$: the correspondence $X \in T_g \mapsto 
\Psi_X \in QD(Y)$ is the Bers embedding $B_Y:T_g \to QD(Y)$.  
As the name suggests, it is an embedding  \cite{Be64}
of the $3g-3$-dimensional Teichm\"uller space $T_g$ into the 
$3g-3$-dimensional complex vector space $QD(Y)$; the 
point $Y$ maps to the origin, and it follows from
results of Nehari \cite{Ne49} that the image is contained in a ball of
radius $6$ and contains a ball of radius $\frac32$.

Within the space $QF$ of quasi-Fuchsian manifolds, the family
$\{Q(X,Y)| X\in T_g\}$ is known as the {\it Bers slice of $QF$
based at Y}.

\sub{2.2 Grafting, Thurston Metric}
Recall that a {\it (complex) projective structure} 
on $S$ is a maximal atlas of
charts from $S$ into $\bbc P^1$ such that all transition maps
are restrictions of elements of $PSL(2,\bbc)$
(i.e. a $(PSL(2,\bbc), \bbc P^1)$-structure in the sense of
Thurston).
Such a structure yields in the usual way a holonomy representation
$hol : \pi_1(S) \to PSL(2,\bbc)$  and an equivariant developing map
$dev : \tilde{S} \to \bbc P^1$.
We will write $P_g$ for the moduli space of projective structures on
$S$ (as defined and topologized, for instance, in \cite{Go88}).

Let $\SS$ denote the set of isotopy classes of essential simple closed
curves on $S$.   There is a well-defined intersection pairing
$i : \SS \times \SS \to \bbz$ given by the minimum number of
intersection points among pairs of representative curves 
in the isotopy classes.  This in turn defines an embedding
of $\bbr_+ \times \SS$ into $\bbr^\SS$ by sending a weighted
simple closed curve $(s,\g)$ to the $\SS$-tuple 
$(s \cdot i(\g, \alpha))_{\alpha \in \SS}$.
The space of measured laminations $\Cal{ML}$ is defined to
be the closure of $\bbr_+ \times \SS$ in $\bbr^\SS$.
For simplicity, a measured lamination coming from a pair $(s, \g)$ will
be denoted $s \g$.

In the presence of a hyperbolic structure on $S$,
it is typical to define measured laminations in terms of
geodesic laminations equipped with a measure on transverse
arcs (see \cite{Th82} or \cite{Bo88} for more details).
We can also use a hyperbolic structure on $S$ to define 
a notion of the length $L(\lambda)$ of a measured lamination
$\lambda$: one defines $L(s\g)$ to be the product of
$s$ and the hyperbolic length of $\g$ on $S$, and then 
extends $L: \bbr_+ \times \SS \to \bbr_+$ to all of 
$\SM \SL$ by continuity (see e.g. \cite{Ke85}).

In \S1, grafting was defined in terms of a map
$\Theta : \Cal{ML} \times T_g \to P_g$; 
for laminations in the 
subset $\bbr_+ \times \SS$ of weighted simple closed curves 
the projective structure $\Theta(s \g, \s)$ was defined
by gluing together the Fuchsian projective structure
associated to $\s$ and a projective $s$-annulus along $\g$.
The proof that $\Theta$ extends continuously to all of 
$\Cal{ML} \times T_g$
can be found in \cite{KT92}.

In order to understand the surjectivity of $\Theta$, let us 
briefly recall the {\it canonical stratification}
associated to a projective structure
(originally due to Thurston -- 
see also \cite{KP94}, \cite{Ap88}, \cite{Sc96}, \cite{KT92}).
First note that, via the developing map $dev$, $\tilde{S}$ inherits
a notion of {\it open round ball} from $\bbc P^1$.
Furthermore, also using $dev$, we can pull back the usual metric
on $\bbc P^1$ to an (incomplete) metric on $\tilde{S}$ --
the metric completion depends only on the projective structure
and is called the {\it M\"obius completion} of $\tilde{S}$ \cite{KP94}.
The closure of an open round ball in the M\"obius completion
is conformally equivalent to compactified hyperbolic space 
$\bbh^2 \cup S^1_\infty$, so the usual notion of 
``hyperbolic convex hull'' transfers.  
Thus, given an open round ball, we write $C(U)$ for the intersection
of $U$ and the convex hull of $\bar{U} \setminus U$ in $\bar{U}$.
The key observation is the following:
\proclaim{Lemma 2.2.1} \cite{KP94}
For every $p \in \tilde{S}$, there is a unique open
round ball $U_p$ such that $p \in C(U_p)$.
\endproclaim
The sets $U_p$ given by the lemma are called 
{\it maximal balls}, and define
a stratification of $\tilde{S}$ into the sets $C(U_p)$
(this descends in turn to a stratification of $S$).
It is easy to verify that in the case of a projective structure
obtained by grafting along $\lambda \in \Cal{ML}$,
this stratification is the basically the same as the one given 
by the leaves and complementary regions of $\lambda$.

We also obtain a
canonical Riemannian metric defined to be the restriction
to $C(U_p)$ of the hyperbolic metric on the
open round ball $U_p$ \cite{KP94}.
We call this metric the {\it grafted metric} or
the {\it Thurston metric}; if the projective structure
is obtained by grafting the hyperbolic surface $\s$
along the measured lamination $\lambda$ (i.e. $\Theta(\lambda,\s)$)
we write this metric as $gr_{\lambda}(\s)$.
Chasing through
the definitions in the case of grafting
along a weighted simple closed curve $s\g$, one sees that
$gr_{\lambda}(\s)$ coincides with $\s$ on $S \setminus \g$ 
and is flat on the inserted annulus.

\sub{2.3 Harmonic maps from surfaces} Let $(M,\s |dz|^2)$ and 
$(N,\rho(w)|dw|^2)$
denote $M$ and $N$ equipped with smooth Riemannian structures; here
$z$ refers to a local conformal coordinate on the surface $M$, and
$w$ refers to a local conformal coordinate on the surface $N$. 
For a Lipschitz map
$w: (M,\s |dz|^2)\to(N,\rho(w)|dw|^2)$, we define the energy 
$E(w ; \s, \rho)$ of
the map $w$ to be
$$
\split
E(w; \s, \rho) &= \int_M \f12 \|dw\|^2 dv(\s)\\
&=\int_M\f1{\s(z)} \{\|w_*\p_z\|^2_{\rho} + \|w_*\p_{\bar z}\|^2_{\rho}\} \s(z)
dzd\bar z
\endsplit
$$
Evidently, while the total energy depends upon the metric structure
of the target surface $(N,\rho)$, it only depends upon the conformal
structure of the source $(M,\s(z))$.

A critical point of this functional is called a harmonic map. We
will be interested in the situation where $M=N=S$, a fixed surface
of finite analytic type, with a fixed homotopy
class $w_0: S\to S$ of maps into the target $S$, where $(S,\rho)$
is possibly singular, but non-positively curved in the 
sense of Alexandrov.  
In that case, (see \cite{GS92;Lemma~1.1, Theorem~2.3}) there
is a unique (if $w_*(\pi_1M)$ is non-abelian) harmonic map $w(\s):
(S,\s)\to(S,\rho)$ in the homotopy class of $w_0$; in the next
section, we will specialize to a case where we will find additional
smoothness for $w$.

For harmonic maps $w: (\SR,\s)\to(N,\rho)$ from a Riemann surface $\SR$
to a smooth target, one can characterize the harmonicity of $w$ in
terms of conformal objects on $\SR$. The pullback metric $w^*{\rho}$
decomposes by type as
$$
\align
w^*{\rho}  &= \left<w_*\p_z, w_*\p_z\right>_{\rho} dz^2 
+ (\|w_*\p_z\|^2_{\rho} +
\|w_*\p_{\bar z}\|^2_{\rho})\s dzd\bar z + \left<w_*\p_{\bar z}, w_*
\p_{\bar z}\right>_{\rho} d\bar z^2\\
&= \vp dz^2 + \s e(w)dzd\bar z + \bar\vp d\bar z^2
\endalign
$$
where $e=\f12\|dw\|^2$ is the energy density of the map $w$. It is
easy to show (see \cite{Sa78}) that if $w$ is harmonic then $\Phi=\vp
dz^2$ is a holomorphic quadratic differential on $\SR$. In particular,
Schoen \cite{Sc84} has emphasized that even for harmonic maps to 
singular metric spaces $(S, \rho)$, it is a consequence
of Weyl's lemma that the differential
$\Phi=\vp dz^2 =  \left<w_*\p_z, w_*\p_z\right>_{\rho} dz^2$ is
holomorphic.

The expression $\SH=\|w_*\p_z\|^2_{\rho}$ plays a special role
in harmonic maps between surfaces (see, for instance \cite{Wo91a}). 
First, we can rewrite the pullback metric
$w^*{\rho}$ entirely in terms of $\Phi=\vp dz^2$ and $\SH$ as follows:

$$
w^*{\rho} = \vp dz^2 + (\SH + \frac {|\Phi|^2}{{\s} {\SH}})dzd\bar z +
\bar {\vp} d\bar z^2.
$$

Moreover,  the function $\SH = \SH(z)$ satisfies the Bochner equation
(this is basically a Liouville equation for prescribed curvature,
using the harmonic map gauge)

$$
\Delta_{\s}\log\SH(z) = -2K_{\rho}(w(z))\{\SH(z) - 
\frac {|\Phi(z)|^2}{{\s(z)} {\SH(z)}}\}
+ 2K_{\s}(z). \tag 2.3.1
$$

Here $K_{\rho}$ and $K_{\s}$ refer to the Gauss curvatures of 
$(S, \s)$ and $(S,\rho)$, respectively, and we are stating the
equation only in the context of smooth maps; we will later 
extend the meaning of this equation to the singular context 
which is our principal interest in this paper.

\sub{2.3.2 Smoothness of Harmonic Maps Families} We will be
interested in harmonic maps between surfaces equipped with the
grafted (Thurston) metrics; in particular, we will carefully study
one-parameter families of such maps. This study relies on the
background result that these maps are reasonably smooth, and that
the family of maps is reasonably smooth in the family parameter, for
a smooth family of grafted metrics. In this section, we establish
these basic smoothness results: the proofs are completely
straightforward generalizations of those found in the literature
(see \cite{Jo97}, \cite{EL81}, \cite{Sa78}), but as the precise versions we
need do not seem to be present already in print, we include them
here for the sake of completeness.

First let us record the regularity of the Thurston metrics;
a proof can be found in \cite{KP94}.

\proclaim{Lemma 2.3.1} For $\la\in\SM\SL$, the grafted metrics 
$\gr_\la(\s)$ are of class $C^{1,1}$.
\endproclaim

Next, we consider the regularity of an individual harmonic map
$w:(S,\gro)\to(S,\gr(\s_1))$.

\proclaim{Lemma 2.3.2} There exists a harmonic map
$w:(S,\gro)\to(S,\gr(\s_1))$ homotopic to the identity; this map
is of class $C^{2,\a}$.
\endproclaim

\prf  As $S$ is compact, and $\gr(\s_1)$ is an NPC space (see
\cite{GS92}), it is straightforward that there is an energy minimizer
$w$ in the given homotopy class. Then we are able to make
considerable use of the literature: Theorem~2.3 of \cite{GS92} then
ensures that $w\in H^1(S,S)$ is locally Lipschitz. The rest of the
proof is straightforward bootstrapping applied to the harmonic map
equation (see, e.g. \cite{Jo97}, 
proof of Theorem~3.2.4).
\qed

Finally, we come to the smoothness of the families of the maps. We
begin by recording the 
fact that $\gr_\la(\s_t)$ varies analytically in $t$, for an analytic
family of hyperbolic metrics $\s_t$.

\proclaim{Lemma 2.3.3} \cite{Mc98} Let
$\{\s_t\}$ be an analytic family (in $t$) of hyperbolic metrics.
Then the family $\{\gr_\la(\s_t)\}$ of grafted metrics is also
analytic in $t$. 
\endproclaim

We omit the proof. Consider such an analytic family
$\{\gr_\la(\s_t)\}$ and the family $\{w_t\}$ of harmonic maps
$w_t:(S,\grt)\to(S,\gro)$ which we know to exist and be of class
$C^{2,\a}$.

\proclaim{Lemma 2.3.4} The family $w_t:(S,\grt)\to(S,\gro)$ of
harmonic maps is analytic in $t$, for small values of $t$.
Any individual map $w_t:(S,\grt)\to(S,\gro)$ is a homeomorphism.
\endproclaim

\prf  We mimic an allied proof in \cite{EL81}: see also \cite{Sa78}.
Given such a family, the first variation at $t=0$ of the tension
$\tau_t=\tau(w_t)$ can be computed to be
$$
\frac{d\tau}{dt} = \Delta\dot w + K_0\dot w -
\text{Tr}\(\frac{\p\G(t)}{\p t}\).\tag2.3.2
$$
$\G(t)$ refers to the Christoffel symbols of the family $\grt$ and
where we have simplified the formula considerably by applying it at
$t=0$, where $w_0:(S,\gro)\to(S,\gro)$ is the identity map.

We aim to apply the analytic
implicit function theorem (see \cite{Be77}): the formal setting is
that we regard the tension $\tau$ as a functional
$$
\tau:C^{2,\a}(S,S)\x(-\e,\e)\lra C^{0,\a}(T(S))
$$
where $C^{0,\a}(T(S))$ denotes $C^{0,\a}$ sections of the tangent
bundle to $S$, the map associates to a map $w\in C^{2,\a}(S,S)$
and a metric $\grt$ the tension field $\tau(w,\grt)$ of the map
$w:(M,\grt)\to(M,\gro)$. This functional is evidently analytic in
$t$, so our attention turns to formula (2.3.2): we assert that
$$
\Big\|\frac{d\tau}{dt}\Big\| > 0
$$
where the norm is that taken on functionals between $C^{2,\a}(S,S)$
and $C^{0,\a}(T(S))$. It is enough to prove that $(\Delta+K_0)$ is
invertible on $C^{0,\a}(T(S))$; i.e. that given $f\in C^{0,\a}(T(S))$,
there is a $u\in C^{2,\a}(S,S)$ so that $\Delta u+K_0u=f$. As
$K_0\le0$, this result follows from standard estimates: e.g.
\cite{GT83, Theorems~8.3, 8.8} give estimates on $\|u\|_{W^{2,2}}$ in
terms of $\|f\|_{L^2}$, and since $\dim_\BR S=2$, this yields a
$C^\a$ estimate on $u$, with higher regularity following from
bootstrapping as in Lemma~2.3.2.

That an individual map $w_t:(S,\grt)\to(S,\gro)$ is a homeomorphism
follows from the map $w_t$ being a perturbation of the identity.
\qed

\rem Here we restrict to families of grafted metrics where the
grafting locus $\lambda$ remains fixed in $t$. If we were to vary
the grafting locus $\lambda= \lambda_t$ in $\SM \SL$, we would
need to deal with issues arising from  $\SM \SL$ having but
a piecewise linear structure and not a differentiable structure.

\sub{2.4 Variation of geodesics}

This section contains a brief discussion of the equations
governing the variation fields of a geodesic in family of
conformally related Riemannian metrics.   
We begin by setting some notation. Consider a smooth family of 
Riemannian metrics 
$g_t$ on $S$ and a family of $g_t$-geodesics $\g_t:[0,1] \to S$. 
We adopt Fermi coordinates
along the curve $\g_0$ so that
$$
g_0=F(x_2)^2dx^2_1+dx^2_2.
$$
The geodesic
equation for $\g_t$ in these coordinates is given by
$$
\g^k_{t,11} + \G^k_{t,ij}(\g_t(x_1))\g^i_{t,1}\g^j_{t,1} = 0\tag2.4.1
$$
where $\G^k_{t,ij}$ are the $g_t$-Christoffel symbols. We
differentiate (2.4.1) in time $t$ to obtain the following equation for
the vector field $\dot\g^k\p_k=\frac d{dt}\bigm|_{t=0}\g^k_t\p_k$
$$
\split
\dot\g^k_{11} + &\(\frac d{dt}\Bigm|_{t=0}\G^k_{t,ij}(\g_0(x_1))\)
\g^i_{0,1}\g^j_{0,1} + D_m\G^k_{0,ij}(\g_0(x_1))\dot\g^m
\g^i_{0,1}\g^j_{0,1}\\
&+\G^k_{0,ij}(\g_0(x_1))\dot\g^i_1\g^j_{0,1} + \G^k_{0,ij}(\g_0(x_1))
\g^i_{0,1}\dot\g^j_1 = 0.
\endsplit
$$
In the Fermi coordinates chosen, we have that $\G^k_{0,ij}=0$ and for a
constant speed geodesic, we have $\g^i_{0,1}=\ell\d^i_1$ where $\ell$ is
the length of the geodesic and $\d^i_1$ is the Kronecker delta. 
Thus the previous equation simplifies to
$$
\dot\g^k_{11} + \(\frac d{dt}\Bigm|_{t=0}\G^k_{t,11}\)\ell^2 +
\(D_m\G^k_{0,11}\)\dot\g^m\ell^2 = 0.\tag2.4.2
$$
We are principally interested in the normal component of the
variation field $\frac d{dt}\g_t$, so we set $k=2$ and compute
$$
\G^2_{t,11} = \frac12g^{2\a}_t\(2\p_1g_{t,1\a} - \p_\a g_{t,11}\).
$$
Moreover, we will be interested only in the situation where
$g_t=\frac1{\SH_t}g_0$ is a family of conformal metrics (see \S3.1) 
and where $g_0$, being written in
Fermi coordinates, is diagonal; this also forces $g_t$ to be diagonal
which simplifies the above description to 
$$
\G^2_{t,11} = -\frac12\SH_t\(\p_2((F(x_2)^2)/\SH_t)\).
$$
It is then straightforward to compute from this equation and from
$\SH_0\equiv1$ that
$$
D_m\G^2_{11} = K\d^2_m\tag2.4.3
$$
and
$$
\aligned
\frac d{dt}\Bigm|_{t=0}\G^2_{t,11}  &= -\frac12\dot\SH
\(\p_2((F(x_2)^2)/\SH_0\)\\
&\qquad-\frac12\lb\{\[-\p_2((F(x_2))^2)\dot\SH/\SH^2_0 + 2(F(x_2))^2
(\p_2\SH_0)\dot\SH/\SH^3_0\]\rb\}\\
&\qquad + \frac12(F(x_2))^2
\frac{\p_2\dot\SH}{\SH^2_0}\\
&= \frac12\p_2\dot\SH
\endaligned\tag2.4.4
$$
where the first and second terms vanish because $\p_2 F(x_2)=0$, and
the third term vanishes because $\p_2\SH_0=\p_2(1)=0$.

We conclude from (2.4.2), (2.4.3) and (2.4.4) that
the variational field $V=\frac d{dt}\bigm|_{t=0}\g^2_t$ 
satisfies
$$
V_{11} + K_0 V\ell^2 = -\frac12\ell^2\p_2\dot\SH.\tag2.4.5
$$
\rem The reader should recognize how, in the case of the Thurston metric
defined above where $K_0$ is
a discontinuous function, the equation (2.4.5) is really a pair of
equations for a single variational field $V$. That is, on the flat
cylinder $K_0\equiv0$ while on the hyperbolic portion of the surface
$K_0\equiv-1$; in our solution for $\dot\SH$ below, this is reflected in 
a jump in the normal derivative of $\dot\SH$ across the 
two geodesics bounding the grafted cylinder.

\sub{\S3. The Case of Simple Closed Curves}

In \S3, we prove the main theorem in the model case when
the measured lamination is a weighted simple closed curve.
We begin by describing the problem in terms of harmonic
maps and deriving our basic equation of study (3.1.5).
The proof effectively becomes a computation,
which we undertake in \S3.2.
As noted earlier, our setup applies quite generally to families of
grafted metrics in which the length of the inserted annulus 
is allowed to vary.   We only use the information that
this length is constant in $t$ at the very end of the proof -- 
this is the content of section \S3.3.

We begin with a precise statement of our objective.

\proclaim{Theorem 3.1} {\rm (Model Case)}. Let $S$ be a closed
differentiable surface of genus $g>1$, let $\g$ be an essential simple closed
curve on $S$ and let $s\in\BR_+$ be a positive real number. Then the
grafting map $\grg:T_g\to T_g$ is a homeomorphism.
\endproclaim

\pf{Theorem 3.1} 

As discussed in \S 1.2, we need only show

\proclaim{Lemma 3.2} The grafting map $\grg$ is locally injective.
\endproclaim

\pf{Lemma 3.2} We suppose, in order to obtain a contradiction, that
there is a hyperbolic surface $(S,\s)$ and an approximating sequence
$\<(S,\s_n)\>\to (S,\s)$ of hyperbolic surfaces so that
$\grg(\s_n)=\grg(\s)$ for all $n$. By passing to a subsequence, and
using the 
differentiability of $\grg:T_g\to T_g$, we find a
tangent vector $[\dot\mu]\in T_{[\s]} T_g$ so that
$d\grg(\s)[\dot\mu]=[0]\in T_{[\grg(\s)]}T_g$.

With this in mind, it is psychologically convenient to solve instead
a formally easier problem: we imagine a differentiable family
$\<(S,\s_t)\>$ of hyperbolic surfaces converging to
$(S,\s)=(S,\s_0)$ with the property that the tangent vector to the
family is given by
$$
\frac d{dt}\Bigm|_{t=0}\<(S,\s_t)\> = [\dot\mu]\in T_{\s_0}T_g
$$
and that $\grg(\s_t)=\grg(\s_0)$ in $T_g$.
We then seek a contradiction to this situation.

Our method is to use harmonic maps to ``fix the gauge'' in comparing
the surfaces $\grg(\s_t)$. In particular, we imagine $\grg(\s_t)$ as
being realized by a metric $\grt$ on the underlying differentiable
surface $S$. Of course, we have a choice for these representative
metrics, as the group $\diffo$ of diffeomorphisms isotopic to the identity
acts on metrics on $S$, with the orbit of $\grt$ consisting of
isometric metrics. However, by the results in \S2, and because
$\gro$ is non-positively curved, there is a unique harmonic map
$w:(S,\grt)\to(S,\gro)$ homotopic to the identity
for any of our choices of $\grt$. In
particular, we can choose this representative metric $\grt$ on $S$
so that the identity map
$$
\id: \(S,\grt\)\lra\(S,\gro\)
$$
is harmonic for all $t\ge0$. Since conformal maps are always
harmonic and our harmonic map is unique, we may conclude
that the identity map above is both harmonic and conformal.

Let
$$
\SH_t = \|\id_{z_t}\|^2 = \frac{\gr(\s_0)}{\grt}|\id_{z_t}|^2\tag3.1.1
$$
denote the holomorphic energy density of the harmonic conformal map,
where here we have snuck in the local coordinate convention that the
metric $\grt$ admits an expansion in the local conformal coordinates
$z_t$ (themselves smooth in $t$) as $\grt=\grt|dz_t|^2$. Then
because the identity map is conformal, we conclude that
$$
\grt|dz_t|^2 = \frac1{\SH_t}\gro|dz_0|^2.\tag3.1.2
$$
Furthermore, because the identity map is harmonic, we apply the
Bochner equation (2.3.1) to conclude that
$$
\Delta_{\grt}\log\SH_t = -2K_0(\id(z))\SH_t + 2K_t(z).\tag3.1.3
$$
We then use (3.1.2) to rewrite (3.1.3) as
$$
\SH_t\Delta_{\gro}\log\SH_t = -2K_0(z)\SH_t + 2K_t(z).
$$
We then divide by $\SH_t$ to obtain the equation
$$
\Delta_{\gro}\log\SH_t = -2K_0 + \frac{2K_t}{\SH_t}\tag3.1.4
$$
which is the precursor
to our basic equation of study. To obtain the basic equation of
study, we differentiate equation (3.1.4) in time to obtain an equation
for $\dot\SH=\frac d{dt}\bigm|_{t=0}\SH_t$
$$
\Delta_{\gro}\dot\SH/\SH_0 = 2\dot K/\SH_0 - 2K_0\dot\SH/\SH^2_0
$$
where $\dot K$ denotes the measure $\frac d{dt}\bigm|_{t=0}K_t$ (see
extended discussion below). Since $\SH_0\equiv1$ by construction, we
summarize our equation as
$$
\(\Delta_{\gro} + 2K_0\)\dot\SH = 2\dot K.\tag3.1.5
$$
This equation requires some discussion. The term $2\dot K$ refers to
a measure supported on the pair of images of the geodesic $\g$ which
bound the grafted cylinder. We can imagine $\dot K$ being
constructed as follows. Since we have a well-defined family of
metrics $\grt$ each with a pair $\{\g^\ell_0,\g^r_0\}$ (left and
right) boundaries of the grafted cylinder, we see that there is then
a pair of well-defined families of curves $\{\g^\ell_t,\g^r_t\}$ of
left and right boundaries of the grafted cylinders. These curves a
priori may vary smoothly on the surface $S$, and form the frontier
of the support of the function $K_t$ which is identically $-1$ on
the hyperbolic portion of $\grt$. Of course, on any open set $\SO$
which avoids the family $\g_t$ for some interval of time
$t\in[0,\e)$, we have that $\frac{d}{dt} K_t=0$ on $\SO$. Thus, $\dot K$ is
supported only on the pair $\{\g^\ell_0,\g^r_0\}$. It is
straightforward, but not required for our work here, to compute
$\dot K$ in terms of the variational vector fields $\{\frac
d{dt}\bigm|_{t=0}\g^\ell_t,\frac d{dt}\bigm|_{t=0}\g^r_t\}$ along
the curves $\{\g^\ell_0,\g^r_0\}$.

In fact, our solution $\frac d{dt}\SH_t$ to equation (3.1.5) is
analytic away from $\{\g^\ell_0,\g^r_0\}$ and Lipschitz in a
neighborhood of those curves. To see that, begin with equation (3.1.1)
and observe that, by definition, the harmonic map
$\id:(S,\grt)\to(S,\gro)$ is conformal hence (locally) complex
analytic (in the local coordinates). Thus, the expression $\frac
d{dt}\id_{z_t}$ is real analytic for an analytic path of metrics
$\{\s_t\}$. This fact implies that when we take
$\ddt\SH_t=\ddt\frac{\gro}{\grt}|\id_{z_t}|^2$, the only possibly
non-analytic contribution comes from the term $\frac
d{dt}\grt(z_t)=\frac\p{\p z}\gro\cd\frac d{dt}z_t$. Finally, from
Lemma~2.3.1, we see that $\frac\p{\p z}\gro$ is Lipschitz continuous
in a neighborhood of $\{\g^\ell_0,\g^r_0\}$, and therefore, so is
$\dot\SH$.
Thus, much
depends on understanding the ``jump'' in derivatives of a solution
$\dot\SH$ across the curves $\{\g^\ell_0,\g^r_0\}$.

\sub{3.2 Computation of $\dot \SH$} 

The goal of this section is a proof of
\proclaim{Lemma 3.3} Any function $\frac d{dt}\bigm|_{t=0}\SH_t$
which solves (3.1.5) and infinitesimally solves (3.1.2) must vanish
identically on $S$. Thus $\id:(S,\grt)\to(S,\gro)$ is an isometry,
up to order $\bigo(t^2)$.
\endproclaim

Recall the basic plan from \S1.2:
We consider equation (3.1.5) as really a pair of equations, the
first saying that $\dot\SH$ is harmonic on the Euclidean cylinder,
and the second saying that $\dot\SH$ satisfies the linearized
Liouville equation on the hyperbolic portion of $\gro$. We then
solve for the general expression for a harmonic function on the
cylinder, and this determines both boundary values for $\dot\SH$ and
normal derivatives $\p_n\dot\SH$ on the cylinder. These derivatives
$\p_n\dot\SH$ on the cylinder can be used to find general solutions
$V^\ell$ and $V^r$
to (2.4.5) with $K\equiv0$ which is compatible with our general
solution to (3.1.5) on the cylinder, and then the version of (2.4.5)
with $K\equiv-1$ gives a general expression for $\p_n\dot\SH$ as
viewed from the hyperbolic side of $\gro$.

We then take this general solution to (3.1.5) and integrate
$\dot\SH\Delta_{\gro}\dot\SH$ by parts to find that equation (3.1.5)
forces $\dot\SH$ to vanish identically. This will prove Lemma~3.3.

We now carry out this outline. To begin, write equation (3.1.5) as
$$
\Delta_E\dot\SH = 0\tag"(3.1.5)$_0$"
$$
and
$$
(\Delta_h - 2)\dot\SH = 0\tag"(3.1.5)$_{-1}$"
$$
where $\Delta_E$ and $\Delta_h$ denote the Laplace-Beltrami
operators on the Euclidean $(K\equiv0)$ and hyperbolic $(K\equiv-1)$
open submanifolds of $(S,\gro)$, respectively. Let us solve
(3.1.5)$_0$ by writing the Euclidean grafted cylinder as
$\{(x,y)\mid-\dcyl \le x\le \dcyl, 0\le y\le\ell\}$ and then writing
$$
\dot\SH = \Sigma a_n(x)\iny.
$$
Then equation (3.1.5)$_0$ implies that
$$
\Sigma\(a''_n(x) - \frac{4\pi^2n^2}{\ell^2}a_n(x)\)\iny = 0.
$$
We conclude that
$$
a_n(x) = \cases
c_n\cosh\frac{2\pi n}{\ell}x + d_n\sinh\frac{2\pi n}{\ell}x  &n\neq0\\
c_0 x + d_0    &n=0
\endcases
$$
so that, for $z=x+iy$,
$$
\dot\SH(z) = c_0 x + d_0 + \Sigma' \(c_n\cosh\frac{2\pi n}{\ell}x +
d_n\sinh\frac{2\pi n}{\ell}x\)\iny.\tag3.2.1
$$
Since $\dot\SH(z)$ is real, we find that
$$
\aligned
\barcn  &= c_{-n}\\
\bardn  &= - d_{-n}.
\endaligned\tag3.2.2
$$
When we specialize (3.2.1) to the boundaries $x=\pm \dcyl$ of the
Euclidean cylinder, we obtain the expansions
$$
\aligned
\dot\SH\Bigm|_{x=-\dcyl}  &= -c_0\dcyl + d_0 + \Sigma'
\(c_n\copi - d_n\sipi\)\iny\\
\dot\SH\Bigm|_{x=\dcyl}  &= c_0\dcyl + d_0 + \Sigma'
\(c_n\copi + d_n\sipi\)\iny.
\endaligned\tag3.2.3
$$
As remarked in the last paragraph of \S3.1,
$\dot \SH$ is real analytic away from the boundary of the
inserted cylinder
and Lipschitz across the boundary.   In particular, the 
$x$-derivative of our general solution (from the cylinder side) 
exists; we compute it by differentiating (3.2.1):
$$
\frac{\p\dot\SH}{\p x}=
\dot\SH_x = c_0 + \Sigma'(\frac{2\pi n}{\ell}c_n\sinh\frac{2\pi
n}{\ell}x + \frac{2\pi n}{\ell}d_n\cosh\frac{2\pi n}{\ell}x)\iny
$$
which we specialize to $x=\pm \dcyl$ to obtain
$$
\aligned
\dot\SH_x\Bigm|_{x=-\dcyl}  &= c_0 + \Sigma'\frac{2\pi n}{\ell}
\(-c_n\sipi + d_n\copi\)\iny\\
\dot\SH_x\Bigm|_{x=\dcyl}  &= c_0 + \Sigma'\frac{2\pi n}{\ell}
\(c_n\sipi + d_n\copi\)\iny.
\endaligned\tag3.2.4
$$
Corresponding to the pair (3.1.5)$_0$ and (3.1.5)$_{-1}$ of versions of
(3.1.5) there is a pair of versions of equation (2.4.5); we intend to
rewrite those equations in terms of the $xy$-coordinates, which 
requires some interpretation beforehand. First we decide
that the $x_2$-direction in \S2.4 will be intepreted as the
$\f\p{\p x}$ direction and the $x_1$-direction will be intepreted as
the $\f\p{\p y}$ direction in a neighborhood of our cylinder; this
is permissible despite forcing $\{\f\p{\p x_1},\f\p{\p x_2}\}$ to be a
frame with the opposite orientation than $\{\f\p{\p x},\f\p{\p y}\}$.
(All this choice of direction amounts to in (2.4.5) is a choice of
direction normal to $\g^\ell$ and $\g^r$ in which $V$ is measured
positively.) Next we observe that the arclength parameter $x_1$
was defined on the domain
$[0,1]$, while the coordinate $y$, which we are presently using to
parametrize the geodesics $\g^{\ell}$ and $\g^r$, 
varies over the domain $[0,\ell]$; we conclude that 
$V_{yy}= \ell^{-2}V_{11}$.
Thus we translate equation (2.4.5) to:
  
$$
\align
V_{yy}  &= -\frac12(\p_x\dot\SH)_0\tag"(2.4.5)$_0$"\\
V_{yy} -V  &= -\frac12(\p_x\dot\SH)_{-1}.
\tag"(2.4.5)$_{-1}$"
\endalign
$$
where here we have written the pair of derivatives of $\dot\SH$ as
$(\p_x\dot\SH)_0$ and $(\p_x\dot\SH)_{-1}$ depending on which side
of $\g^\ell$ or $\g^r$ we are considering. (In this notation,
formula (3.2.4) refers to $(\p_x\dot\SH)_0$.)  Once again, we are
using the regularity of the solution $\dot \SH$ for the existence
of the $x$-derivatives
$(\p_x\dot\SH)_0$ and $(\p_x\dot\SH)_{-1}$.

At this point, we need to consider that there are two boundary
components $\g^\ell$ and $\g^r$, and hence two variation vector
fields with normal components $V_-$ defined along
$\g^\ell=\{x=-\dcyl\}$ and $V_+$ defined along $\g^r=\{x=+\dcyl\}$. If
we then set
$$
V_- = \Sigma\la_n\iny\quad\text{and}\quad
V_+ = \Sigma\rho_n\iny\tag3.2.5
$$
(where $\la_n$ and $\rho_n$ suggest ``left'' and ``right'',
respectively), we substitute (3.2.4) and (3.2.5) into (2.4.5)$_0$ to obtain
$$
\Sigma\frac{-4\pi^2n^2}{\ell^2}\la_n\iny = \frac{-c_0}2 - \frac12
\Sigma'\frac{2\pi n}\ell\(-c_n\sipi + d_n\copi\)\iny.
$$
We solve for $\la_n$ in the above to obtain, for $n\neq0$,
$$
\la_n = \frac{\ell}{4\pi n}\(-c_n\sipi + d_n\copi\).
\tag"(3.2.6)$_-$"
$$
Similarly, we find, for $n\neq0$,
$$
\rho_n = \frac{\ell}{4\pi n}\(c_n\sipi + d_n\copi\).
\tag"(3.2.6)$_+$"
$$
Along the way, we also find that, setting $n=0$,
$$
c_0 = 0.\tag3.2.7
$$
The formulae (3.2.6)$_-$ and (3.2.6)$_+$ when substituted into (3.2.5)
yield the expansions
$$
\aligned
V_-  &= \la_0 + \Sigma'\frac{\ell}{4\pi n}\(-c_n\sipi + d_n\copi\)\iny \\
\text{and}\quad
V_+  &= \rho_0 + \Sigma'\frac{\ell}{4\pi n}\(c_n\sipi + d_n\copi\)\iny.
\endaligned
$$
We now use the crucial observation that $V_-$ (and $V_+$, of course)
solves both equations (2.4.5)$_0$ and (2.4.5)$_{-1}$. Thus, from our
knowledge of $V_-$ and $V_+$, we can apply (2.4.5)$_{-1}$
to this expansion and obtain $(\p_x\dot\SH)_{-1}$. That is,
(2.4.5)$_{-1}$ is equivalent to
$$
\aligned
&(\p_x\dot\SH)_{-1}\Bigm|_{x=-\dcyl} = -2\{(V_-)_{yy} -
(V_-)\}\\
&= 2\la_0 -2\Sigma'\frac{\ell}{4\pi n}\(-c_n\sipi + d_n\copi\)
\(-\frac{4\pi^2n^2}{\ell^2} - 1\)\iny \\
&= 2\la_0 +\Sigma'\frac1{2\pi n\ell}\(4\pi^2n^2 + \ell^2\)
\(-c_n\sipi + d_n\copi\)\iny.
\endaligned\tag"(3.2.8)$_-$"
$$
Similarly
$$
(\p_x\dot\SH)_{-1}\Bigm|_{x=+\dcyl} = 
2\rho_0 +
\Sigma'\frac1{2\pi n \ell}
\(4\pi^2n^2 + \ell^2\)\(c_n\sipi + d_n\copi\)\iny.\tag"(3.2.8)$_+$"
$$
We pause and observe that we have obtained in equations (3.2.3)
and (3.2.8) the Dirichlet and Neumann conditions, respectively, for
the linearized Liouville/Bochner equation (3.1.5)$_{-1}$. This permits
us to focus for the rest of the computation on the compact
hyperbolic surface $(S_{-1},\gro|_{S_{-1}})$ where $S_{-1}$ denotes
the closure of the $\{K_0=-1\}$ subdomain of $S$.

It is a reflex in this situation to integrate by parts:
$$
\aligned
0  &= \iint\limits_{S_{-1}}\dot\SH\Delta_h\dot\SH -
2\dot\SH^2dA_{\gro}\\
&= \iint\limits_{S_{-1}} -\|\nabla\dot\SH\|^2 - 2\dot\SH^2dA_{\gro} +
\int\limits_{\p S_{-1}}\dot\SH(\p_n\dot\SH)ds_{\gro}
\endaligned\tag3.2.9
$$
where $\p_n\dot\SH$ denotes the outward normal. The boundary term we
compute with (3.2.3) and (3.2.8), using that the normal 
derivatives
$\p_n\dot\SH\bigm|_{x=-\dcyl}=\p_x\dot\SH$ and
$\p_n\dot\SH\bigm|_{x=+\dcyl}=-\p_x\dot\SH$
in our coordinates.  Thus
$$
\aligned
&\int\limits_{\p S_{-1}}\dot\SH(\p_n\dot\SH)ds_{\gro} = 
\int\limits_{\{x=-\dcyl\}}\dot\SH(\p_x\dot\SH)dy +
\int\limits_{\{x=\dcyl\}}\dot\SH(-\p_x\dot\SH)dy\\
&= 2\ell d_0\la_0 + \int\limits^{y=\ell}_{\Sb y=0\\ \{x=-\dcyl\}\endSb}\Sigma'
\frac1{2\pi n \ell}\(4\pi^2n^2 + \ell^2\)\(-c_n\sipi + d_n\copi\)\times\\
&\qquad\(c_{-n}\cosh\frac{\pi(-n)s}{\ell} -
d_{-n}\sinh\frac{\pi(-n)s}{\ell}\) +\\
&\qquad\{\text{terms involving non-trivial powers of $\iny$}\} dy +\\
&\qquad-2\ell d_0\rho_0 + 
\int\limits^{y=\ell}_{\Sb y=0\\ \{x=+\dcyl\}\endSb}\Sigma'(-1)
\frac1{2\pi n \ell}\(4\pi^2n^2 + \ell^2\) \(c_n\sipi + d_n\copi\)\times\\
&\qquad\(c_{-n}\cosh\frac{\pi(-n)s}{\ell} +
d_{-n}\sinh\frac{\pi(-n)s}{\ell}\) +\\
&\qquad\{\text{terms involving non-trivial powers of $\iny$}\} dy
\endaligned
$$
after substituting (3.2.3) and (3.2.8). After integration and applying
(3.2.2), the above expression simplifies to
$$
\aligned
&\int\limits_{\p S_{-1}}\dot\SH(\p_n\dot\SH)ds_{\gro} = 
2\ell d_0(\la_0-\rho_0) +\\
&\ell\Sigma'
\frac1{2\pi n \ell}\(4\pi^2n^2 + \ell^2\)\(-c_n\sipi + d_n\copi\)
\(\barcn\copi - \bardn\sipi\)\\
&- \ell\Sigma'\frac1{2\pi n \ell}\(4\pi^2n^2 + \ell^2\)\(c_n\sipi + d_n\copi\)
\(\barcn\copi + \bardn\sipi\)\\
&= 2\ell d_0(\la_0-\rho_0) +\\
&\frac1{2\pi}\Sigma'\frac1n\(4\pi^2n^2 + \ell^2\)
\(-(|c_n|^2+|d_n|^2)\sipi\copi + c_n\bardn \sinh^2\frac{\pi ns}{\ell} +
\barcn d_n \cosh^2\frac{\pi ns}{\ell}\)\\
& -\frac1{2\pi}\Sigma'\frac1n\(4\pi^2n^2 + \ell^2\)
\((|c_n|^2+|d_n|^2)\sipi\copi + c_n\bardn \sinh^2\frac{\pi ns}{\ell} +
\barcn d_n \cosh^2\frac{\pi ns}{\ell}\)\\
&= 2\ell d_0(\la_0-\rho_0)
-\frac1{\pi} \Sigma'\frac1n\(4\pi^2n^2 + \ell^2\)
(|c_n|^2 + |d_n|^2)\sipi\copi.
\endaligned\tag3.2.10
$$
Thus this final integral is the negative of the sum of positive 
terms summed with a mystery term $2\ell d_0(\la_0-\rho_0)$.
In Lemma~3.4, we will use the ``slice condition'' that our family of
grafted metrics have Euclidean cylinders of unvarying length $s$
to conclude that this term $2\ell d_0(\la_0-\rho_0)$ is non-positive,
which will force the integral in (3.2.10) to be non-positive.

\proclaim{Lemma 3.4} $\la_0 -\rho_0 = -\frac{sd_0}{2}$
\endproclaim

We postpone the proof of Lemma~3.4 until \S3.3, preferring to 
assume it for now to finish the proof of Lemma~3.3. 
Now, assuming Lemma~3.4 and applying (3.2.10) to (3.2.9), we
find that

$$
\iint\limits_{S_{-1}} - \|\nabla\dot\SH\|^2 - 2\dot\SH^2dA_{\gro} -
\frac1{\pi}\Sigma'\frac1n\(4\pi^2n^2 + \ell^2\)(|c_n|^2 + |d_n|^2)\sipi\copi 
-{\ell}sd_0^2= 0\tag"(3.2.11)$^*$"
$$
with all terms being nonpositive: this forces $c_n=d_n=0$ for all $n$ 
and $\dot\SH$ to vanish
identically in $S_{-1}$.  Because $\dot\SH$ is continuous
across the boundary of the cylinder, it vanishes on all of $S$.
Thus $\grt=\gro(1+\bigo(t^2))$, so that
$\id:(S,\grt)\to(S,\gro))$ is an isometry, up to order $\bigo(t^2)$.
This concludes the proof of Lemma~3.3.
\qed

From (3.2.11)$^*$ we also obtain,
as a corollary of the proof,

\proclaim{Lemma 3.5} For a family of conformal grafted metrics
$\grt$, we have $V_+\equiv V_-\equiv0$.
\endproclaim

\prf We see that $c_n=d_n=0$, so formulae (3.2.6) show that
$\lambda_n = \rho_n = 0$ for $n \ne 0$. But then also equation
(2.4.5)$_{-1}$ forces $\lambda_0 = \rho_0 =0$.
\qed

It remains to conclude the proof of Lemma~3.2. We already
have from Lemma~3.3 that $\grt=\gro(1+\bigo(t^2))$; thus
$\id:(S,\grt)\to(S,\gro)$ is an isometry, up to order $\bigo(t^2)$.

Now, if $(S,g)$ is a metric space which can be written as the graft
of a hyperbolic surface with geodesic boundary $\g$ to a flat cylinder of
height $s$, then there is a well-defined inverse operation
$\ungrg$ to the grafting operation $\grg$: the inverse
operation excises the flat cylinder and reidentifies the hyperbolic
surface with boundary along $\g$ by projecting the
flat cylinder along its longitudinal geodesics. Of course, we have the
property that $\ungrg \circ \grg$ is the identity isometry.
We apply this operation $\ungrg$ to $(S,\grt)$: because
$\id:(S,\grt)\to(S,\gro)$ is an isometry to order $\bigo(t^2)$,
and Lemma~3.5 guarantees that the geodesic $\g$ is not moving
(to $\bigo(t^2)$),
we find that the hyperbolic portions $(S_{-1}, \gro)$ and
$(S_{-1}, \grt)$ agree to $\bigo(t^2)$.  Moreover, the flat
portions $(S_{0}, \gro)$ and
$(S_{0}, \grt)$ also agree to $\bigo(t^2)$, and so again 
citing Lemma~3.5, we conclude that the longitudinal geodesics
across the flat cylinders are also unchanged, to $\bigo(t^2)$.
Thus $\s_0 = \ungrg\{\grt + \bigo(t^2)\} = 
\ungrg\{\grt\}  + \bigo(t^2) = \s_t  + \bigo(t^2)$, proving
the lemma.
\qed

\sub{3.3 Slice Condition} We have yet to use the hypothesis that the
grafted cylinder has constant length $s$ in the family $\grt$.
Certainly it is necessary to use this hypothesis to prove Lemma~3.2,
as \tec space is $6g-6$ (real) dimensional and the space
$T_g\x\BR_+$ of grafted hyperbolic metrics (up to $\diffo$
equivalence) is $6g-5$ (real) dimensional. Thus we might expect that
the map $T_g \x \BR_+ \to T_g$ which records the conformal equivalence class
of an equivalence class of grafted metrics would pullback points to
one-dimensional families of grafted metrics. The content to
Lemma~3.2 is that such families would meet level sets
$T_g \times \{s_0\}\subset T_g \times \BR_+$ in points; thus we must
somehow use the fact that we are restricted to such a level set in
the proof of Lemma~3.2.

Let us extend the notation of \S3.1 somewhat and allow the Euclidean
portion of the grafted metric $\grt$ to be a Euclidean cylinder of
length $s=s(t)$, where we permit the length $s(t)$ to vary in $t$.
In the notation of \S3.2, we claim that

\proclaim{Lemma 3.6}
$\la_0-\rho_0=-\frac{sd_0}{2}-\frac{ds(t)}{dt}$.
\endproclaim

From this lemma, we see the

\pf{Lemma 3.4} The Euclidean length $s(t)$ is constant, hence
$\frac{ds(t)}{dt}=0$, and result follows.

\pf{Lemma 3.6} The plan is to compute the $t$-derivative of the area
$A(\grt)$ of the family $\grt$ of grafted metrics two ways. In the
first method, we use that $\grt$ is a metric which is composed of a
portion which is hyperbolic with geodesic boundary and a portion
which is composed of a Euclidean cylinder, and so the area is
compatible via Gauss-Bonnet and elementary geometry. In the second
method, we use the analytical formulae (3.1.2) and (3.1.5).

\sub{First Method} If we remove the cylindrical portion of the
grafted metric and glue the resulting hyperbolic
surface-with-geodesic-boundary together across its pair of geodesic
boundary components, we obtain a closed hyperbolic surface of area
$-2\pi\chi$, where $\chi$ is the Euler characteristic. Thus, using
that the cylinder has length $\elg\cd s(t)$, we find that the area
$A(\grt)$ of the grafted metric satisfies
$$
A(\grt) = - 2\pi\chi + \elg\cd s(t)
$$
so that
$$
\ddt A(\grt) = \(\ddt\elg\)s(0) + \elg\cd\ddt s(t).\tag3.3.1
$$
To find the derivative $\ddt\elg$, we first observe that since the
length $\elg$ is that of a geodesic $\g_t$ which varies smoothly in
a family containing the {\it geodesic} $\g_0$, then we must have that
$$
\ddt\elg = \ddt\int_{\g_0}\dst.\tag3.3.2
$$
Yet the term $\dst$ is computable from (3.1.2) as
$$
\dst = \sqrt{\frac1{\SH_t}}\ds.
$$
We combine this equation with (3.3.2) and differentiate to find
$$
\split
\ddt\elg  &= -\frac12\int_\g\dot\SH\SH^{-3/2}_0\ds\\
&= -\frac12\int\dot\SH\ds
\endsplit
$$
as $\SH_0\equiv1$.

We next apply our formulae (3.2.3) for $\dot\SH$ and (3.2.7) for $c_0$
to this last equation to find that
$$
\aligned
\ddt\elg  &= -\frac12\int^\ell_0d_0 +\text{ terms involving powers
of }\exp(2\pi iny/\ell)dy\\
&= -\frac12d_0\ell.
\endaligned\tag3.3.3
$$
Combining (3.3.1) and (3.3.3) yields
$$
\frac d{dt}A(\grt) = -\frac12d_0\ell s(0) + \ell\ddt s(t).\tag3.3.4
$$

\sub{Second Method} Formula (3.1.2) suggests another method, as the
area $A(\grt)$ may be expressed as
$$
\split
A(\grt)  &= \iint\limits_S dA(\grt)\\
&= \iint\limits_S\frac{dA(\gro)}{\SH_t}.
\endsplit
$$

Thus
$$
\split
\ddt A(\grt)  &= -\iint\limits_S\dot\SH\SH^{-2}_0 dA(\gro)\\
&= -\iint\limits_{S_{-1}}\dot\SH dA(\gro)
- \iint\limits_{S_0}\dot\SH dA(\gro)
\endsplit\tag3.3.5
$$

The two terms in formula (3.3.5) require separate treatments. To
evaluate the first term, begin with equation (3.1.5)$_{-1}$ and
integrate to find
$$
0 = \iint\limits_{S_{-1}}(\Delta_h - 2)\dot\SH dA_{\gro} =
\int_{\p S_{-1}}\p_n\dot\SH ds_{\gro} - 2
\iint\limits_{S_{-1}}\dot\SH dA_{\gro}.
$$
We rearrange to find
$$
\split
-\iint\limits_{S_{-1}}\dot\SH  &= -\frac12\int_{\p S_{-1}}
\p_n\dot\SH ds_{\gro}\\
&= -\frac12\int_{\{x=-\dcyl\}}\dot\SH_x dy + \frac12
\int_{\{x=\dcyl\}}\dot\SH_x dy
\endsplit
$$
(by $\p_n$ referring to the outward normal)
$$
\aligned
&= -\frac12\int^\ell_0 2\la_0 +\text{ \{terms involving
non-trivial powers of }\exp(2\pi iny/\ell)\}dy\\
&  +\frac12\int^\ell_0 2\rho_0 +\text{ \{terms involving
non-trivial powers of }\exp(2\pi iny/\ell)\}dy
\endaligned
$$
(by (3.2.8)$_-$ and (3.2.8)$_+$)
$$
= - \ell(\la_0 -\rho_0).\tag3.3.6
$$
To find the second term in (3.3.5), we simply use formula (3.2.1),
again using (3.2.7) to set $c_0=0$. We find
$$
\aligned
-\iint_{S_0}\dot\SH dA_{\gro} &= -\int^{\dcyl}_{-\dcyl}\int^\ell_0 d_0 +
\text{ \{terms involving non-trivial powers of }\exp(2\pi
iny/\ell)\}\\
&= -d_0\ell s(0).
\endaligned\tag3.3.7
$$
We combine (3.3.5), (3.3.6) and (3.3.7) to find
$$
\frac d{dt} A(\grt) = -\ell(\la_0 -\rho_0) - d_0\ell s(0).\tag3.3.8
$$

\sub{Summary} Formulae (3.3.4) and (3.3.8) combine to yield
$$
-\frac12 d_0\ell s(0) + \ell\ddt s(t) = -\ell(\la_0 -\rho_0) -
d_0\ell s(0)\tag3.3.9
$$
from which the statement of the lemma follows immediately.
\qed

\sub{\S4 The General Case}

In this section we will prove the main theorem 
in the case of grafting on a measured lamination 
which is not necessarily a weighted simple closed curve.

\proclaim{Theorem} For any $\lambda \in \Cal{ML}$, 
$Gr_{\lambda} : T_g \to T_g$ is a homeomorphism.
\endproclaim

As in the model case, we need only prove the local injectivity; to
that end, 
we suppose the theorem is false
and get a variation $\s_t$ of $\s_0$ such that
$Gr_{\lambda} \s_t = Gr_{\lambda} \s_0$ in $T_g$.
The harmonic map setup is the same as the model case;
we isotope the grafted metrics and assume that
the identity map $(S, gr_\lambda(\s_t)) \to (S, gr_\lambda(\s_0))$
is harmonic and conformal for each $t$.   The functions
$\SH_t$ and $\dot \SH$ are defined in the usual way.

The first step in the proof, of course, is to approximate 
$\lambda \in \SM\SL$
by a sequence of weighted simple closed curves
$s_m \g_m \to \lambda$ and attempt to use our computations
from the model case.   
The main difficulty is that the family of grafted metrics
$Gr_{s_m \g_m}(\s_t)$ can no longer be assumed conformal
and we must generalize some results from \S2.4 
and \S3 to allow for this possibility.

Our plan is to 
carry out the derivation of \S3 for a single
non-conformal deformation; 
thus we will suppress the subscripts $m$ in the
notation until section \S4.5. 
Let
us observe that there were four basic steps in \S3: 1) The derivation of the
equation (3.1.5) for $\dot\SH$, 2) the computation of $(\p_x\dot\SH)_{-1}$ using
the geodesic variational vector field equation (2.4.5), 3) the relating of the
difference $\la_0-\rho_0$ of the constant terms in the variational vector field
to global quantities in grafting and $\dot\SH$, and 4) the derivation
of the identity for the norm of $\dot\SH$ in terms of $\ell$, $s$ and
$\f{ds}{dt}$, as embodied in formulae (3.2.9), (3.2.10) and the 
formula
(3.2.11)$^*$. We carry out these steps in the next four sections
\S\S4.1-4.4, culminating in a formula like that in formula (3.2.11)$^*$.
We will then interpret this generalized formula
(3.2.11)$^*$ in terms of quantities which are continuous on $\SM\SL$.

\sub{4.1 The Infinitesimal Bochner Equation} In this section,
we show that the basic global equations are unchanged.

\proclaim{Lemma 4.1.1} $\dot \SH$ satisfies equation (3.1.5).
\endproclaim

Thus our main equation of study is unchanged despite allowing
the family of harmonic maps to stray from conformality.
\prf
To see this, we begin with the Bochner equation (2.3.1)
$$
\Delta_t \log \SH_t = -2 K_0 \SH_t + 
2 K_0 \frac{t^2 |\Phi|^2}{\rho_t \SH_t} + 2 K_t
$$
and differentiate once with respect to $t$ at $t = 0$.
We see immediately that 
$$
\frac{d}{dt} 2 K_0 \frac{t^2 |\Phi|^2}{\rho_t \SH_t} = 0
$$
so the derivative of the right hand side becomes
$-2K_0 \dot \SH + 2 \dot K$.
Because $\log \SH_0 \equiv 0$,
all terms of $\frac{d}{dt} \Delta_t$ vanish except those for
which the derivative passes through.
Thus
$$ 
\frac{d}{dt} \Delta_t \log \SH_t = \Delta_0 \dot \SH,
$$
yielding (3.1.5).
\qed

\sub{4.2. The geodesic variational vector field} A straightforward computation
confirms the formula
$$
\grt = -\Phi(t)dz^2 + \frac{\gro}{\SH_t} - \overline{\Phi(t)}d\bar z^2 
+ \bigo(t^2).\tag4.2.1
$$

From this formula, we can compute the expressions we need in order to apply
formula (2.4.1) to the present case. In particular, when we differentiate
(2.4.1) in time, we observe that formula (2.4.2) (and those
formulae following it) 
continues to be valid; we are
left to evaluate $\ddt\G^2_{t,11}$, where we adopt Fermi coordinates for
$\gro$. In terms of those Fermi coordinates $\{x^1,x^2\}$ for $\gro$ along the
curve $\g$ (here representing one of $\{\g^\ell,\g^r\}$), the formula (4.2.1)
becomes
$$
\grt = \frac{\gro}{\SH_t} + t\phi_{ij}dx_idx_j + \bigo(t^2)\tag4.2.2
$$
where the tensor $\phi_{ij}dx_idx_j$ may be represented as
$$
\phi_{ij} = \pmatrix  -2\ree\phi   &2\im\phi\\
2\im\phi   &2\ree\phi
\endpmatrix + \bigo(x^2_2)\bigo(t) + \bigo(t^2).\tag4.2.3
$$
Here we have written 
$\Phi(t)=\{t\phi + \bigo(t^2)\}\(dx_1+i\frac{dx_2}{F(x_2)}\)^2$, noting that
$dx_1+i\frac{dx_2}{F(x_2)}$ is a conformal coordinate up to order 
$\bigo(x^2_2)$. We continue to consider both $\Phi(t)$ and $\phi$
as small quantities, since we regard $gr(\s_t)$ as nearly conformal to
$gr(\s_0)$ to first order in $t$.

Thus, since 
$F(x_2) = 1 + \bigo(x_2^2)$, we may compute along the curve $\g$ that,
$$
\split
\ddt\G^2_{t,11}  &= \frac12\p_2\dot\SH + \frac12\[2\frac\p{\p
x_1}(2\im\phi) + 
\frac\p{\p x_2}(2\ree\phi)\] + \bigo(x^2_2)\\
&= \frac12\p_2\dot\SH + \frac\p{\p x_1}\im\phi
\endsplit
$$
where we evaluate along the curve $\g=\{x_2=0\}$ and we use the Cauchy-Riemann
equations for $\phi$ to simplify the second term. 
Next, jump to section~3.2 where
we use these equations to solve for $(\p_x\dot\SH)_{-1}$ in terms of
$(\p_x\dot\SH)_0$: we adapt our work there to the present (more
general) 
case by
expanding the harmonic function $\im\phi$ in the expansion (compare (3.2.1))
$$
\im\phi = u_0 x + v_0 + \Sigma' \(u_n\cosh\frac{2\pi n}{\ell}x +
v_n\sinh\frac{2\pi n}{\ell}x\)\iny.\tag4.2.4
$$
Let $W$ solve the amended variational equation
$$
W_{yy} = -\frac12(\p_x\dot\SH)_0 - (\p_y\im\phi) \tag"(2.4.5)$^{\Phi}_0$"
$$
in analogy to equation (2.4.5)$_0$.

We find that if $W$ is of the forms $W_-=\Sigma\la^*_n\iny$ and
$W_+=\Sigma\rho^*_n\iny$, then for $n\neq0$
$$
\align
\la^*_n  &= \la_n + \frac{\ell}{2\pi in}\(u_n\copi - v_n\sipi\).
\tag"(4.2.5)$_-$"\\
\rho^*_n  &= \rho_n + \frac{\ell}{2\pi in}\(u_n\copi + v_n\sipi\).
\tag"(4.2.5)$_+$"
\endalign
$$
again finding that
$$
c_0 = 0.
$$

Next, in analogy to (3.2.8), we find that
$$
(\p_x\dot\SH)_{-1} = (\p_x\dot\SH)_0 - 2W\tag4.2.6
$$
and so has the expansions
$$
\aligned
&(\p_x\dot\SH)_{-1}\Bigm|_{x=-\dcyl} = 2\la_0 + \Sigma'(4\pi^2n^2+\ell^2)
\bigg\{\frac1{2\pi n\ell} \(-c_n\sipi + d_n\copi\)\\
&+ \frac1{\pi in\ell} \(u_n\copi - v_n\sipi\)\bigg\}
\endaligned\tag"(4.2.7)$_-$"
$$
and
$$
\aligned
&(\p_x\dot\SH)_{-1}\Bigm|_{x=\dcyl} = 2\rho_0 + \Sigma'(4\pi^2n^2+\ell^2)
\bigg\{\frac1{2\pi n\ell} \(c_n\sipi + d_n\copi\)\\
&+ \frac1{\pi in\ell} \(u_n\copi + v_n\sipi\)\bigg\}
\endaligned\tag"(4.2.7)$_+$"
$$
We will use these formulas in \S4.4 when we combine all of 
our modifications
to \S3 to get a new version of formula (3.2.11)$^*$.

\sub{4.3. The slice condition} In this section, we generalize \S3.3 to
the
(general) case of a non-conformal deformation. Here the metric $\grt$ is 
involved
in the computation of the first variation of arclength (3.3.3) and in the
computation of the first variation of area (3.3.5). 
In the latter, we note from (4.2.1) that
$$
dA(\grt) = \frac{dA(\gro)}{\SH_t(1-|\nu(t)|^2)}
$$
where $\nu(t)=\bigo(t)$ denotes the Beltrami differential. Thus
$$
dA(\grt) = \frac{dA(\gro)}{\SH_t} + \bigo(t)^2.\tag4.3.1
$$
In the computation of the first variation of 
arclength, we have, from (4.2.2) and (4.2.3) that
$$
\split
\ddt\elg  &= \ddt\int_{\g_0}\dst\\
&= \frac d{dt}\int_\g\sqrt{\SH^{-1}_t+2\ree t\phi+\bigo(t^2)}\ ds_{\gro}\\
&= -\frac12\int_\g(\dot\SH - 2\ree\phi)ds_{\gro}\\
&= \frac12\int_\g\dot\SH ds_{\gro} + \ell \bigo(\|\dot\Phi(t)\|).
\endsplit\tag4.3.2
$$
The last expression may require some explanation. For any given infinitesimal
holomorphic quadratic differential $\dot\Phi(t)$, on the fixed compact surface
$\gro$, the Harnack inequality bounds the 
supremum of $|\ree\phi|$ over the curve
$\g$ in terms of the integral norm of $\dot\Phi(t)$ and the length of the curve.
As there is but a compact set of such unit norm quadratic 
differentials $\dot\Phi(t)$, the result holds
in general, even in a precompact family $\{\g\}$ of grafting loci.

Finally we collect terms, as in the summary (3.3.9), 
but with the addition of the
considerations from (4.3.1) and (4.3.2), to find

\proclaim{Lemma 4.3.1} $\ell(\la_0-\rho_0)=-\frac12d_0\ell
s(0)-\ell\frac{ds}{dt}+\ell s(0)\bigo(\|\dot\Phi\|)$.
\endproclaim

\sub{4.4. The extended identity} Finally, we combine 
the results of our previous
sections into an identity analogous to that of (3.2.11)$^*$. 
We extend the computation
of (3.2.9) and (3.2.10) to obtain
$$
\aligned
&\int\limits_{\p S_{-1}}\dot\SH(\p_n\dot\SH)ds_{\gro}  = 
2\ell d_0(\la_0-\rho_0) \\ 
&\qquad - \frac1\pi\Sigma'\frac1n\(4\pi^2n^2 +
\ell^2\) (|c_n|^2+|d_n|^2)\sipi\copi \\
&+ \ell\Sigma'\(4\pi^2n^2 + \ell^2\)\frac{1}{\pi in\ell}\bigg\{\(u_n\copi -
v_n\sipi\) \(\bar c_n\copi - \bar d_n\sipi\)\\
&- \(u_n\copi + v_n\sipi\) \(\bar c_n\copi + \bar d_n\sipi\)\bigg\}\\
&= 2\ell d_0(\la_0-\rho_0) - \frac1\pi\Sigma'\frac1n\(4\pi^2n^2 + \ell^2\)
(|c_n|^2+|d_n|^2)\sipi\copi\\
&\qquad- \Sigma'\(4\pi^2n^2 + \ell^2\)\frac2{\pi in} (v_n\bar c_n + u_n
\bar d_n)
\sipi\copi\\
&= 2\ell d_0(\la_0-\rho_0) - \frac1\pi\Sigma'\frac1n\(4\pi^2n^2 + \ell^2\)
(|c_n|^2+|d_n|^2)\sipi\copi\\
&\qquad- \sum_{n>0}\(4\pi^2n^2 + \ell^2\)\frac2{\pi n}\(\im\[v_n\bar c_n+u_n\bar
d_n\]\)\sipi\copi.
\endaligned
$$
When we combine this with Lemma 4.3.1, we find the
analogue to (3.2.11)$^*$:
$$
\aligned
&\iint\limits_{S_{-1}} - \|\nabla\dot\SH\|^2 - 2\dot\SH^2dA_{\gro} -
\frac1{\pi}\Sigma'\frac1n\(4\pi^2n^2 + \ell^2\)(|c_n|^2 + |d_n|^2)\sipi\copi\\
&- \sum_{n>0}\(4\pi^2n^2 + \ell^2\)\frac2{\pi n}\(\im\[v_n\bar c_n+u_n\bar
d_n\]\)\sipi\copi\\
&- \ell s d^2_0-2\ell sd_0\frac{\dot s}s + \ell s\bigo(\|\dot\Phi\|) = 0.
\endaligned\tag"(4.4.1)$^*$"
$$

\sub{4.5. An identity on $\SM\SL$ and the conclusion of the proof} We define the
length $L=L(\g)$ for the grafted curve $\g$ to be
$$
L(\g) = \ell_\g\cd s;
$$
if we regard the transverse measure on $\g$ as being given by $s\cd i(0,\g)$,
then we see that $L(\g)$ is simply the length of $[\g]$ as an element of the
measured lamination space $\SM\SL$. It is well-known then that the function
$L:\SS\x\BR_+ \subset\SM\SL\to\BR$ extends to a continuous function
$L:\SM\SL\to\BR$. Our main observation is that (4.4.1)$^*$ extends to an identity in
terms of $\dot\SH$, $\dot\Phi$, $\frac d{dt}(\log s)$ and quantities that are
continuous on $\SM\SL$. To see this we rewrite (4.4.1)$^*$ as
$$
\aligned
&\iint\limits_{S_{-1}} - \|\nabla\dot\SH\|^2 - 2\dot\SH^2dA_{\gro} -
\frac1{\pi}\Sigma'\frac1n\(4\pi^2n^2 + \ell^2\)(|c_n|^2 + |d_n|^2)
\sinh\frac{\pi nL}{\ell^2}\cosh\frac{\pi nL}{\ell^2}\\
&- \sum_{n>0}\(4\pi^2n^2 + \ell^2\)\frac2{\pi n}\(\im\[v_n\bar
c_n+u_n\bar d_n\]\)\sinh\frac{\pi nL}{\ell^2}\cosh\frac{\pi nL}{\ell^2} -
Ld^2_0 - 2Ld_0\frac d{dt}\log s\\
&+ L\bigo(\|\dot\Phi(t)\|) = 0.
\endaligned\tag"(4.5.1)$^*$"
$$

It is now straightforward to complete the proof of the main technical result.

\pf{Theorem A} We have already established the result in the case where $\la$ is
a lamination supported on a finite set $\{\g_1,\dots,\g_k\}$ of simple closed
curves. For the general case consider a sequence $\<s_m \g_m\>$ of measured
laminations supported on simple closed curves of lengths $\ell_n=\ell(\g_n)$
which approximates $\la$. 
In our notation, $s_m$ denotes the multiple of the
transverse (intersection) measure for the geodesic $\g_m$, 
so we may express the transverse
measure for $\la$ as $\lim_{m\to\infty}s_mi(\cd,\g_m)$ (see the discussion
in \S2.2).

In analogy with the opening of the proof of Lemma 3.2, we suppose (in order to
obtain a contradiction) that there is a family of surfaces $\s_t$ so that
$$
\text{Gr}_\la(\s_t) = \text{Gr}_\la(\s_0)\tag4.5.2
$$
at least to an order $\bigo(t^2)$ in $T_g$. We then consider the family
$\text{Gr}_{s_m \g_m}(\s_t)$: it is of course no longer necessary that
$\text{Gr}_{s_m \g_m}(\s_t)$ should equal $\text{Gr}_{s_m \g_m}(\s_0)$ to
$\bigo(t^2)$, but the condition (4.5.2) should be asymptotically true in $m$, by
construction. This implies that the
Hopf differential (which we will denote $\Phi_m(t)$)and  which 
measures the quasiconformality between
$\text{Gr}_{s_m \g_m}(\s_t)$ and $\text{Gr}_{s_m \g_m}(\s_0)$)  
should have first variation $\ddt\Phi_m(t)$ in 
$t$ which tends to zero as $m\to\infty$. We write
$$
\lim_{m\to\infty}\ddt\Phi_m(t) = 0
$$
which implies that, in the formula (4.5.1)$^*$, 
we can take $|u_n|$, $|v_n|\to0$,
along with the final term $L\bigo(\|\Phi\|)$,
since $L$ is bounded. With the exception of the term
involving $\ddt\log s(t)$, all of the term involve quantities which are
continuous on $\SM\SL$: this follows immediately, once we observe that
$\ell^2\sinh\frac{\pi nL}{\ell^2}$ is bounded for,
say, $\ell>1$. Of the
exceptional term, we have the lemma

\proclaim{Lemma 4.5.1} $\ddt\log(s(t))=0$.
\endproclaim

\prf
We simply need to interpret the ``slice condition'' for $\lambda$
correctly.  In particular, suppose we allow the transverse measure 
of $\lambda$ to vary as our grafted metrics vary,
i.e. let $\lambda_t = \rho_t \lambda$
with $\rho_t$ differentiable and $\rho_0 = 0$.
The slice condition is given by 
$$
\dot{\rho} = \ddt \rho_t = 0.
$$
Now we can approximate $\lambda_t$ by scaling the
weighted simple closed curves $\{ s_m \g_m \}$ in exactly the same way:
$$
\lambda_t = \lim_{m \to \infty} \rho_t s_m \g_m.
$$
Having done so, we have $s(t) = \rho_t s_m$ in the calculations above,
and $0 = \dot{\rho} = \frac{\dot{s}}{s_m}$.
\qed

The rest of the proof is straightforward. As $m\to\infty$,
we have $\dot\Phi\to0$ and the expression collapses to the previous version
(3.2.11)$^*$. 
In particular, we see from (4.5.1)$^*$ that 
$\dot\SH_m = \bigo(\|\dot\Phi_m(t)\|)$, where we have added subscripts
to $\dot\SH$ and $\Phi(t)$ to emphasize the dependence of these
quantities upon the approximating sequence. Thus, the argument at the 
end of \S3.2 extends to show that $\s_{t,m}$ agrees with $\s_{0,m}$
to $o_m(1) + \bigo(t^2)$ (where $o_m(1)$ goes to zero as
$m$ tends to infinity), and hence that $\s_t$ agrees with $\s_0$
to $\bigo(t^2)$, as required.
This proves the local injectivity, 
as
desired. \qed

\rem We have unnecessarily restricted ourselves to unpunctured
surfaces, primarily for notational simplicity and expositional
cleanliness. The proofs all extend to the punctured case once we
make three observations: i) all of the measured laminations
under consideration avoid a neighborhood of the punctures,
ii) there is a unique harmonic map of finite
energy between surfaces of bounded non-positive curvature and some
negative curvature \cite{Al64} (\cite{Wo91b}), and (iii) the holomorphic energy
function $\SH$ for such a map is bounded in $C^1$ across the
punctures, so no new non-vanishing boundary terms would arise in a
modified formula (3.2.9), or therefore in any of the subsequent 
starred formulae (3.2.11)$^*$, (4.4.1)$^*$, and (4.5.1)$^*$.

\sub{\S5 Applications}

\sub{5.1. Geometric Coordinates on the Bers Slice} Consider
the Bers slice
$$
B_Y = T_g \times \{ Y \} \subset T_g \times T_g \cong QF. 
$$
There is a natural map $\beta : B_Y \to \Cal{ML}$
which assigns to a quasi-Fuchsian group the bending lamination
$\lambda$ on the component of the convex hull boundary facing the
fixed structure $Y$ (continuity of $\beta$ is proved in \cite{KS95}).   
The hyperbolic structure on this
component of the convex hull boundary defines a point $\mu$
in Teichm\"uller space, and the relevant observation is
that 
$Y = \text{Gr}_\lambda(\mu).$
Theorem A shows that the metric $\mu$ is determined by
$Y$ and $\lambda$; therefore since 
the Thurston homeomorphism
$\Theta: \Cal{ML} \times T_g \to P_g$ (described in the
the introduction \S1) is one-one, 
the map $\beta$ is also.
This is a simple way of assigning ``bending coordinates'' 
to $B_Y$.

\proclaim{Corollary 5.1} Let $B_Y$ be a Bers slice with fixed conformal
structure $Y$.  Then the map assigning
the bending lamination on the component of the convex hull boundary
facing $Y$ is a homeomorphism onto its image.
\endproclaim

\sub{5.2 Generalized Bers Slices. Deformation Spaces
of Books of I-bundles}
Any geometrically finite, freely indecomposable 
Kleinian group $G$ has a space $\ST(G)$ of 
quasiconformal deformations parametrized by the product 
$\prod T(S_i)$ of the
\tec spaces $T(S_i)$ of its boundary components 
$\{S_1, \dots, S_K\}$ at infinity \cite{Ma74}. As in the
example of the quasi-Fuchsian groups above, one can define slices
$\SS(t_1, \dots, t_{K-1})$ of these deformation spaces
by simply fixing, say, the conformal structures at infinity of the
first $K-1$ boundary components, and letting the last
conformal structure vary over its \tec space $T(S_K)$.

Let us focus our attention
on a class of geometrically finite three-manifolds
homeomorphic to the interior of a book of I-bundles. The deformation
spaces of these three-manifolds are studied in detail
in \cite{AC96} and are important because of the discovery by 
Anderson and Canary that the closures of
those deformation spaces exhibit previously unexpected
phenomena. The simplest of
these manifolds has the
following description. 
Begin with a solid torus with three disjoint parallel annuli
on the boundary; here we choose the annuli so that their
central curve is homotopic within the solid torus to the
core curve of the solid torus.  Attach, along those annuli,
thickenings $\{H_1,H_2,H_3\}$ of one-holed surfaces
$\{M_1,M_2,M_3\}$ of genera $h_1,h_2$, and $h_3$ (respectively).
The new three-manifold $N$ has boundary surfaces $\{S_1,S_2,S_3\}$
of genera $h_1+h_2, h_2+h_3$ and $h_3+h_1$; indeed, these
bounding surfaces $S_1,S_2$ and $S_3$ are obtained by gluing
$M_i$ to $M_{i+1}$ (with cyclic indexing) along the single
boundary $\g=\p M_i=\p M_{i+1}$. Because all the thickenings of the
surfaces are glued along neighborhoods of curves which retract to the
core curve, we see that all of the curves $\g_i$ are homotopic to 
each other and to the core curve $\g$ of the central solid torus.

Now consider the space $\ST(G)$ of quasi-conformal deformations 
of the Kleinian group $G$ obtained as the holonomy of the 
hyperbolization of this three-manifold $N$. We consider the
slice $\SS(t_1,t_2) \subset \ST(G) = T(S_1) \times T(S_2) \times
T(S_3)$ defined by the coordinate description 
$$
\SS(t_1,t_2) = \{t_1\} \times \{t_2\} \times T(S_3)
$$
We then consider a map $\b_{AC}:\SS(t_1,t_2) \to \SM\SL \times \SM\SL$,
analogous to the map $\b:B_Y \to \SM\SL$ above, which assigns to an
element $(t_1, t_2, t_3) \in \SS(t_1,t_2)$ the pair of bending 
measures of the boundary components 
($C_1$ and $C_2$)of the convex hulls facing the
conformal structures at infinity represented by $t_1$ and $t_2$.
Our application of Theorem~A is the following

\proclaim{Corollary 5.2} 
The map $\b_{AC}:\SS(t_1,t_2) \to \SM\SL \times \SM\SL$ is injective.
\endproclaim

\prf 
Suppose $\b_{AC}(t_1,t_2,t_3) = \b_{AC}(t_1,t_2,t_3')$. Then
by Theorem~A, not only do the bending measures on the 
convex hull boundary components $C_i$ and $C_i'$ ($i=1,2$,
facing the ends $t_i$ and $t_i'$, respectively), but so do the
hyperbolic structures. Lift to the quasi-Fuchsian covers 
$\SQ_i$ and $\SQ_i'$ of $(t_1,t_2,t_3)$ and $(t_1,t_2,t_3')$
corresponding to the surfaces $C_i$ and $C_i'$ and observe that
these are identical by Corollary~5.1.  Thus the holonomy 
representations of $\pi_1(C_i)$ are conjugate to those
of $\pi_1(C_i')$. But as there is a common element $\g$ in
$\pi_1(C_1) \subset \pi_1(N)$ and $\pi_1(C_2)\subset \pi_1(N)$,
we see that the pair of representations of 
$\pi_1(S_1) *_{[\g]} \pi_1(M_3) = \pi_1(N(t_1,t_2,t_3))$
(in the obvious notation) are conjugate. Thus
$[(t_1,t_2,t_3)]=[(t_1,t_2,t_3')] \in \ST(G)$, as desired.
\qed

\sub{5.3. 2+1 de Sitter Spacetimes} 
We finish with an 
application to the structure of 
(2+1)-dimensional de~Sitter spacetimes, following \cite{Sc96}.  
Recall that $3$-dimensional de~Sitter space is defined to be the set of unit
spacelike vectors in Minkowski space:
$$
\desitter = \{ v \in \bbr^4_1 | < v,v > = +1 \}.
$$
This is the model space for Lorentzian 3-manifolds of constant
positive curvature.   Projectivizing $\bbr^4_1$ to $\bbr P^3$,
we get the Klein model of hyperbolic space from the unit timelike
vectors, the sphere at infinity $S^2_\infty$ from the light cone, and a
projective model of $\desitter$ as the remainder of $\bbr P^n$.
Taking polar duals with respect to the sphere at infinity
gives a correspondence between points in the projectivized de~Sitter space
and planes in hyperbolic space 
(and thus with round circles on $S^2_\infty$).

Now imagine a projective structure on a closed hyperbolic surface $S$
close to a Fuchsian structure (the construction works for any
projective structure but is easiest to describe for the quasi-Fuchsian case).
Using the polarity mentioned above, the set of all closed round balls
contained within $dev(\tilde{S})$ defines a certain open subset
$\Cal{U}$ of $\desitter$.  The holonomy $hol(\pi_1(S))$ acts
discontinuously on $\Cal{U}$ and the quotient is a de~Sitter spacetime
homeomorphic to $S \times \bbr$.   Any example arising from
a projective structure on $S$ in this way is called a 
{\it standard de~Sitter spacetime}.
Standard de~Sitter spacetimes are well-behaved from
the point of view of causality -- in particular, we can choose
the product structure so that each slice $S \times \{t\}$ is
spacelike and every timelike or lightlike curve crosses 
$S \times \{t\}$ exactly once (we say $S \times \bbr$ is
a {\it domain of dependence}). 
  
The main result of \cite{Sc96} is that every de~Sitter spacetime
$S \times \bbr$
which is a domain of dependence embeds in a standard
de~Sitter spacetime.
Now suppose we have an example coming from
a projective structure with Thurston coordinates 
$(\lambda, \sigma) \in \Cal{ML} \times T_g$.
A domain of dependence has a well-defined {\it causal horizon};
it follows easily that the causal horizon corresponds
to the space of {\it maximal} open round balls, which is in turn isometric
to the $\bbr$-tree dual to $\lambda$ \cite{Sc96}.

We are now able to refine our classification
of de~Sitter spacetimes, by providing
coordinates in terms of naturally-arising data in the future
(the future causal horizon) and the past (the conformal structure
on $S$ at past infinity).
More precisely, we have the following reworking of Theorem A:

\proclaim{Corollary 5.3} 
Let $\lambda \in \Cal{ML}$ be a measured lamination with dual
$\bbr$-tree $\tilde{\lambda}$.  Let $\desitter(S ; \tilde{\lambda})$
be the family of standard de~Sitter spacetimes with future causal
horizon $\tilde{\lambda}$, and define a map
$c_\infty : \desitter(S ; \tilde{\lambda}) \to T_g$ which 
assigns the conformal structure on $S$ at past infinity.
Then $c_\infty$ is one-one.
\endproclaim

\prf
By definition, any two standard de~Sitter spacetimes $\Cal{M}_i$ 
($i = 1,2$) in $\desitter(S; \tilde{\lambda})$  come from
projective structures on $S$ (say with Thurston
coordinates $(\lambda, \sigma_i)$).    
By examining the construction above,
we have 
$$
c_\infty(\Cal{M}_i) = Gr_{\lambda}(\sigma_i).
$$
Because $Gr_{\lambda}$ is one-one (Theorem A), $c_\infty$ is also.
\qed

\end